\documentclass[11pt,leqno]{amsart}
\usepackage{amsmath, amssymb, amsaddr}
\usepackage{graphicx,color,hyperref}
\usepackage{verbatim,enumitem}

\def\smallskip{\vskip3pt plus1pt minus1pt}
\def\medskip{\vskip5pt plus2pt minus1pt}

\def\plainsubsection#1|{%
  \par\vskip0.25cm\penalty -100
  \centerline{{\sc #1}}
  \vskip3pt plus 1pt minus 0pt
  \penalty 500}

\long\def\claim#1|#2\endclaim{\par\vskip 4pt\noindent 
{\bf #1.}\ {\em #2}\par\vskip 5pt}

\def\plainproof{\noindent{\em Proof}}

\def\today{\ifcase\month\or
January\or February\or March\or April\or May\or June\or July\or August\or
September\or October\or November\or December\fi \space\number\day,
\number\year}

\catcode`\@=11
\newcount\@tempcnta \newcount\@tempcntb 
\def\timeofday{{%
\@tempcnta=\time \divide\@tempcnta by 60 \@tempcntb=\@tempcnta
\multiply\@tempcntb by -60 \advance\@tempcntb by \time
\ifnum\@tempcntb > 9 \number\@tempcnta:\number\@tempcntb
  \else\number\@tempcnta:0\number\@tempcntb\fi}}
\catcode`\@=12

\def\Bibitem#1&#2&#3&#4&%
{\hangindent=2cm\hangafter=1
\noindent\rlap{\hbox{\bf #1}}\kern2cm{\rm #2}{\it #3}{\rm #4.}} 

\def\bC{{\mathbb C}}

\def\bN{{\mathbb N}}
\def\bP{{\mathbb P}}
\def\bQ{{\mathbb Q}}
\def\bR{{\mathbb R}}

\def\bOne{{\mathchoice {\rm 1\mskip-4mu l} {\rm 1\mskip-4mu l}
{\rm 1\mskip-4.5mu l} {\rm 1\mskip-5mu l}}}

\def\cE{{\mathcal E}}
\def\cF{{\mathcal F}}
\def\cG{{\mathcal G}}
\def\cI{{\mathcal I}}
\def\cJ{{\mathcal J}}
\def\cL{{\mathcal L}}

\def\cO{{\mathcal O}}

\def\cV{{\mathcal V}}

\def\cX{{\mathcal X}}

\def\gm{{\frak m}}

\def\ii{i}
\def\ld{,\ldots,}
\def\bu{{\scriptstyle\bullet}}
\def\ort{\mathop{\hbox{\kern1pt\vrule width4pt height0.4pt depth0pt
    \vrule width0.4pt height7pt depth0pt\kern3pt}}}

\def\qedsquare{\hbox{
\vrule height 1.5ex  width 0.1ex  depth 0ex\kern-0.1ex
\vrule height 1.5ex  width 1.5ex  depth -1.4ex\kern-1.5ex
\vrule height 0.1ex  width 1.5ex  depth 0ex\kern-0.1ex
\vrule height 1.5ex  width 0.1ex  depth 0ex}\kern0.5pt}
\def\qed{~\hfill\qedsquare\vskip6pt plus2pt minus1pt}
\def\lambdawedge{\mathop{\raise1.5pt\hbox{$\scriptstyle\bigwedge$}}\nolimits}

\let\ssm\smallsetminus

\let\le\leqslant
\let\leq\leqslant
\let\compact\Subset
\let\ge\geqslant
\let\geq\geqslant

\def\swt#1{\smash{\widetilde#1}}
\def\swh#1{\smash{\widehat#1}}
\def\build#1^#2_#3{\mathrel{\mathop{\null#1}\limits^{#2}_{#3}}}

\def\mertorelbar{\vrule width0.6ex height0.65ex depth-0.55ex}
\def\merto{\mathrel{\mertorelbar\kern1.3pt\mertorelbar\kern1.3pt\mertorelbar
    \kern1.3pt\mertorelbar\kern-1ex\raise0.28ex\hbox{${\scriptscriptstyle>}$}}}
\let\lra=\longrightarrow
\let\lra=\longrightarrow

\def\lraww{\mathrel{\rlap{$\longrightarrow$}\kern-1pt\longrightarrow}}

\catcode`\@=11
\newdimen\@rrowlength \@rrowlength=6ex
\def\ssrelbar{\vrule width\@rrowlength height0.64ex depth-0.56ex\kern-4pt}
\def\llra#1{\@rrowlength=#1\ssrelbar\rightarrow}
\def\vlra#1{\hbox to#1mm{\rightarrowfill}}
\catcode`\@=12

\def\lcm{\mathop{\rm lcm}\nolimits}

\def\End{\mathop{\rm End}\nolimits}

\def\Tr{\mathop{\rm Tr}\nolimits}

\def\Proj{\mathop{\rm Proj}\nolimits}
\def\Supp{\mathop{\rm Supp}\nolimits}

\def\codim{\mathop{\rm codim}\nolimits}
\def\rank{\mathop{\rm rank}\nolimits}

\def\Gr{\mathop{\rm Gr}\nolimits}

\def\ddbar{{\partial\overline\partial}}
\def\bddK{{{}^b\kern-1pt K}}

\def\Sing{\mathop{\rm Sing}}

\def\tors{{\rm tors}}
\def\reg{{\rm reg}}

\def\Sing{{\rm Sing}}
\def\FS{{\rm FS}}

\def\GG{{\rm GG}}
\def\abs{{\rm abs}}

\newdimen\plainitemindent \plainitemindent=18pt
\def\plainitem#1{\vskip3pt\noindent
\hangindent\plainitemindent\hbox to\plainitemindent{#1\hss}\ignorespaces}

\catcode`\@=11
\def\openup{\afterassignment\@penup\dimen@=}
\def\@penup{\advance\lineskip\dimen@
  \advance\baselineskip\dimen@
  \advance\lineskiplimit\dimen@}
\newdimen\jot \jot=3pt
\newskip\plaincentering \plaincentering=0pt plus 1000pt minus 1000pt
\def\ialign{\everycr{}\tabskip\z@skip\halign}
\def\eqalign#1{\null\,\vcenter{\openup\jot\m@th
  \ialign{\strut\hfil$\displaystyle{##}$&$\displaystyle{{}##}$\hfil
      \crcr#1\crcr}}\,}
\newif\ifdt@p
\def\displ@y{\global\dt@ptrue\openup\jot\m@th
  \everycr{\noalign{\ifdt@p \global\dt@pfalse \ifdim\prevdepth>-1000\p@
      \vskip-\lineskiplimit \vskip\normallineskiplimit \fi
      \else \penalty\interdisplaylinepenalty \fi}}}
\def\@lign{\tabskip\z@skip\everycr{}} 
\def\displaylines#1{\displ@y \tabskip\z@skip
  \halign{\hbox to\displaywidth{$\@lign\hfil\displaystyle##\hfil$}\crcr
    #1\crcr}}
\def\eqalignno#1{\displ@y \tabskip\plaincentering
  \halign to\displaywidth{\hfil$\@lign\displaystyle{##}$\tabskip\z@skip
    &$\@lign\displaystyle{{}##}$\hfil\tabskip\plaincentering
    &\llap{$\@lign##$}\tabskip\z@skip\crcr
    #1\crcr}}
\def\leqalignno#1{\displ@y \tabskip\plaincentering
  \halign to\displaywidth{\hfil$\@lign\displaystyle{##}$\tabskip\z@skip
    &$\@lign\displaystyle{{}##}$\hfil\tabskip\plaincentering
    &\kern-\displaywidth\rlap{$\@lign##$}\tabskip\displaywidth\crcr
    #1\crcr}}
\def\plaincases#1{\left\{\,\vcenter{\normalbaselines\m@th
    \ialign{$##\hfil$&\quad##\hfil\crcr#1\crcr}}\right.}
\def\plainmatrix#1{\null\,\vcenter{\normalbaselines\m@th
    \ialign{\hfil$##$\hfil&&\quad\hfil$##$\hfil\crcr
      \mathstrut\crcr\noalign{\kern-\baselineskip}
      #1\crcr\mathstrut\crcr\noalign{\kern-\baselineskip}}}\,}

\newcommand\myleaders{\leavevmode\leaders\hbox{\kern1pt.}\hfill}
\renewcommand{\tocsection}[3]{%
  \indentlabel{\@ifnotempty{#2}{\ignorespaces#1 
      \rlap{\kern0.1em #2.}\hphantom{00.\ \ }}}#3\myleaders}
\catcode`\@=12

\title[Morse cohomology estimates for jet differential 
operators]{Morse cohomology estimates\\ for jet differential operators}

\author[Jean-Pierre Demailly, Mohammad Reza Rahmati]
{Jean-Pierre Demailly${}^*$, Mohammad Reza Rahmati${}^{**}$}
\address{${}^*$Institut Fourier, Universit\'e Grenoble Alpes, Gi\`eres, France\\
Institute of Algebraic Geometry, Leibniz Universit\"at Hannover, Germany\\
\vskip0.75cm
{\it Dedicated to the memory of Professor Paolo de Bartolomeis}
\vskip0.25cm\strut\vskip0.5cm
}

\date{}
\thanks{Both authors are supported by the ERC project ALKAGE,
contract 670846 from September 2015.}

\begin{document}

\begin{abstract}
We provide detailed holomorphic Morse estimates for the cohomology
of sheaves of jet differentials and their dual sheaves. These
estimates apply on arbitrary directed varieties, and a special
attention has been given to the analysis of the singular
situation. As a consequence, we obtain existence results for global
jet differentials and global differential operators under positivity
conditions for the canonical or anticanonical sheaf of the directed
structure.
\vskip10pt
\noindent
{\bf Key words:} Directed variety, jet bundle, jet differential, jet metric,
holomorphic Morse inequalities, canonical sheaf.
\vskip10pt
\noindent
{\bf MSC~Classification (2010):} 32H30, 32L10, 14J17, 14J40, 53C55
\end{abstract}

\maketitle

\tableofcontents

\setcounter{section}{-1}
\section{Introduction}

Recent developments in the theory of Kobayashi hyperbolicity have
shown that a very convenient framework is the category of directed
varieties.  By definition, a directed variety is a pair $(X,V)$ where
$X$ is a complex manifold or variety, and $V\subset T_X$ a complex analytic
linear subspace that may itself possess singularities. Such a structure
is defined alternatively as a saturated coherent subsheaf
$\cV$ of $\cO(T_X)$.  One can then
introduce a natural concept of canonical sheaf sequence
$\smash{K_V^{[\bullet]}}$ that generalizes the usual canonical sheaf
$\det(V^*)$ in the case where $V$ is nonsingular; one says that
$(X,V)$ is of general type if $\smash{K_V^{[\bullet]}}$ is big (see section~2
for details).

Following classical ideas of A.~Bloch [Blo26a, 26b], Green-Griffiths [GG79]
greatly advanced the study of Kobayashi hyperbolicity through a powerful
use of jet bundles; as an application, they obtained a new geometric
proof of the Bloch conjecture. In this vein, to any directed variety
$(X,V)$ one can associate a bundle $J^kV$ of $k$-jets of holomorphic curves
$f:(\bC,0)\to X$ tangent to $V$. This gives rise to tautological
rank~$1$ sheaves $\cO_{X_k^\GG}(m)$ on the weighted projectivized bundles
$$
\pi_k:X_k^\GG\to X,\qquad X_k^\GG:=(J^kV\ssm\{0\})/\bC^*.
\leqno(0.1)
$$
By taking direct images (cf.\ [GG79] and [Dem95], see also section~3),
one gets sheaves of jet differentials
$$
\cO(E_{k,m}V^*)=(\pi_k)_*\cO_{X_k^\GG}(m).\leqno(0.2)
$$
The goal of this  work is to obtain precise estimates for the dimensions
of the cohomology groups
$$
H^q\big(X,\cO(E_{k,m}V^*)\big)\quad\hbox{and}\quad
H^q\big(X,\cO((E_{k,m}V^*)^*)\big),\leqno(0.3)
$$
along the lines of [Dem11]. The crucial technical argument is an
application of holomorphic Morse inequalities to the rank one sheaf
$\smash{\cO_{X_k^\GG}(m)}$, suitably twisted. The main contribution of
the present exposition is to analyze the role of singularities in a
more systematic way than in our previous papers [Dem11,~12, 15, 18].
The method to cope with singularities is to introduce ad hoc sheaves of
holomorphic sections that are bounded near singular points. In a quite
general setting, these sheaves interact very well with holomorphic Morse
inequalities, and allow us to extend the estimates of the nonsingular case
in a simple manner. Another new idea is to introduce the
dual sheaf $\smash{\cO((E_{k,m}V^*)^*)}$, which can be seen as a sheaf of
differential operators acting on functions of $J^kV$. In this
direction, we prove an existence theorem for twisted differential operators.
It somehow extends the existence theorem of slanted vector fields on 
jet spaces, as established by Siu [Siu04], P\u{a}un [Pau08] and 
Merker [Mer09]. The positivity conditions needed
to get the existence of twisted differential operators are much weaker than
those needed for the existence of twisted vector fields, but
the draw-back of the present approach is that we get a priori no
information on the degeneration loci.
\medskip

Paolo de Bartolomeis was one of the great world experts of deformation
theory.  In the present context, it would be interesting to
investigate the deformation theory of directed varieties, in the
smooth and singular contexts as well. For instance, for a deformation
$(\cX,\cV)\to S$ of directed structures over a base~$S$, a~natural
question is the invariance of ``directed plurigenera''
$\smash{h^0(X_t,{}^bK_{V_t}^m)}$ (where 
$\smash{{}^bK_{V_t}^m}$ refers to the sheaf of
bounded sections, see \S$\,$2), or at least the invariance of 
the volume of $K_{V_t}$ along the fibers $X_t$ of $\cX\to S$.
\medskip
The first author wishes to thank Mihai P\u{a}un for raising several important
issues that play a very significant role in this work.

\medskip

\section{Category of directed varieties}

We first recall the main definitions concerning the category of
directed varieties. We start with the nonsingular case and then explain
in detail the additional concepts and requirements that we introduce in
the presence of singularities.

\claim 1.1. Definition|A $($complex, nonsingular$\,)$ {\rm directed variety} 
is a pair $(X,V)$ consisting of a $n$-dimensional complex manifold $X$ 
equipped with a holomorphic subbundle $V\subset T_X$. 
A morphism \hbox{$\Phi:(X,V)\to (Y,W)$} in the category of directed varieties 
is a holomorphic map such that \hbox{$d\Phi_x(V_x)\subset W_{\Phi(x)}$}
for every point $x\in X$.
\endclaim

The absolute situation is the case $V=T_X$
and the relative situation is the case when $V=T_{X/S}$ is the
relative tangent space to a smooth holomorphic map $X\to S$. In
general, we can associate to $V$ a sheaf $\cV=\cO(V)\subset\cO(T_X)$
of holomorphic sections. No assumption need be made on the Lie bracket tensor
$[\bu,\bu]:\cV\times \cV\to \cO(T_X)/\cV$, i.e.\ we do not assume any
sort of integrability for~$\cV\,$; if~this happens, then $\cV$ defines
a holomorphic foliation on~$X$. 

\claim 1.2. Complement to the definition (singular case)|{\rm Usually we
are interested in questions that are birationally invariant; in such cases
one can always blow-up $X$ and reduce the situation to the case when $X$
is nonsingular. Even then, the tangent bundle $T_X$ need not have
many holomorphic subbundles, but it always possesses many analytic
subsheaves, thus it is important to allow singularities for~$V$. When
defining a directed structure $V$ on $X$, we assume that there exists
a dense Zariski open set
$X'=X\ssm Y\subset X$ such that $V_{|X'}$ is a subbundle of $T_{X|X'}$ and
the closure $\overline{V_{|X'}}$ in the total space of $T_X$ is an analytic
subset. The rank $r\in\{0,1,\ldots,n\}$ of $V$ is by definition the
dimension of $V_x$ at points $x\in X'\,$; the dimension may be larger
at points $x\in Y$. This happens e.g.\ on $X=\bC^n$ for the rank~$1$
linear space $V$ generated by the Euler vector field
$\varepsilon(z)=\bC\sum_{1\le j\le n} z_j{\partial\over\partial z_j}\,$:
then $V_z=\bC\varepsilon(z)$ for $z\ne 0$, and $V_0=\bC^n$.
The singular set $\Sing(V)$ is by definition 
the complement of the largest open subset on which $V$ is a subbundle 
of $T_X\,$; it is equal to the indeterminacy set of the meromorphic
map $X\merto\Gr(T_X,r)$ into the Grassmannian bundle of $r$-dimensional
subspaces of $T_X$, hence $\codim(\Sing(V))\geq 2$. The category of
directed varieties $(X,V)$ is obtained by allowing $X$ to be an arbitrary
reduced complex space; if $X\hookrightarrow Z$ is a local embedding
in a smooth ambient variety $Z$, we assume that there exists a
Zariski open set $X''\subset X'=X_{\reg}$ on which $V_{|X''}$ is 
a subbundle of $T_{X'}$, and that $\overline{V_{|X''}}$ is a closed
analytic subset of~$T_Z$ (similarly $T_X$ is defined to be the closure 
of $T_{X'}$ in $T_Z$).}
\endclaim 

\section{Pluricanonical sheaves of a directed variety}

Let $(X,V)$ be a directed projective manifold where $V$ is possibly singular,
and let $r=\rank V$. If $\mu:\widehat X\to X$ is a proper modification
(a composition of blow-ups with smooth centers, say), we get a directed
manifold $(\widehat X,\widehat V)$ by taking $\widehat V$ 
to be the closure of $\mu_*^{-1}(V')$, where $V'=V_{|X'}$ is the
restriction of $V$ over a Zariski open set $X'\subset X\smallsetminus \Sing(V)$ 
such that $\mu:\mu^{-1}(X')\to X'$ is a biholomorphism. We say that
$(\widehat X,\widehat V)$ is a {\em modification} of
$(X,V)$ and write $\swh{V}=\mu^*V$.

We will be interested in taking modifications realized by iterated blow-ups of 
certain nonsingular subvarieties of the singular set $\Sing(V)$, so 
as to eventually ``improve'' the singularities of $V\,$; outside of
$\Sing(V)$ the effect of blowing-up is irrelevant. The 
canonical sheaf $K_V$, resp.\ the pluricanonical sheaf sequence  $K^{[m]}_V$, 
is defined by using the concept of bounded pluricanonical forms 
that was already introduced in [Dem11].

\claim 2.1. Definition|For a directed pair $(X,V)$ with $X$ nonsingular,
we define $\bddK_V$, resp.\ $\bddK^{[m]}_V$, for any integer $m\ge 0$, 
to be the rank~$1$ analytic sheaves such that
$$
\eqalign{
\bddK_V(U)&=\hbox{sheaf of locally bounded sections 
of}~~\cO_X\big(\Lambda^r V^{\prime *}\big)(U\cap X')\cr
\bddK^{[m]}_V(U)&=\hbox{sheaf of locally bounded sections 
of}~~\cO_X\big((\Lambda^r V^{\prime *})^{\otimes m}\big)(U\cap X')\cr}
$$
where $r=\hbox{rank}(V)$, $X'=X\smallsetminus \Sing(V)$, $V'=V_{|X'}$,
and ``locally bounded'' means bounded with respect to a smooth hermitian 
metric $h_X$ on $T_X$, on every set $U_c\cap X'$ such that $U_c$ is relatively
compact in $U$.
\endclaim

The above definition of $\bddK_V^{[m]}$ may look like an
analytic one, but it can easily be turned into an equivalent algebraic
definition (cf.\ [Dem18]). Let us recall that, given a coherent ideal sheaf
$\cJ=(g_1,\ldots,g_N)$ and a positive rational (or even real) number $p$,
one defines the $p$-th integral closure, denoted formally $\overline{\cJ^p}$,
to be the sheaf of holomorphic functions $f$ that satisfy locally an inequality
$|f|\leq C\big(\sum|g_j|)^p$; this is a coherent sheaf that can be
identified with the sheaf of functions satisfying an integral
equation $f^d+a_1f^{d-1}+\ldots+a_d=0$ where $a_s\in\cJ^{\lceil ps\rceil}$
for some $d\geq 1$ (see e.g.\ [agbook]).

\claim 2.2. Proposition|Consider the natural morphism 
$\cO(\Lambda^rT_X^*)\to \cO(\Lambda^r V^*)$ where $r=\rank V$ and
$\cO(\Lambda^r V^*)$ is defined as the quotient of 
$\cO(\Lambda^rT_X^*)$ by $r$-forms that have zero restrictions 
to $\cO(\Lambda^rV^*)$ on $X\smallsetminus \Sing(V)$. The 
bidual $\cL_V=\cO(\Lambda^r V^*)^{**}$ is an invertible sheaf,
and our natural morphism can be written
$$
\cO(\Lambda^rT_X^*)\to \cO(\Lambda^rV^*)=\cL_V\otimes\cJ_V\subset \cL_V
\leqno(2.2\,{\rm a})
$$
where $\cJ_V$ is a certain ideal sheaf of $\cO_X$ whose zero set is
contained in $\Sing(V)$ and the arrow on the left is surjective by 
definition. Then 
$$
\bddK_V^{[m]}=\cL_V^{\otimes m}\otimes\overline{\cJ^m_V}\leqno(2.2\,{\rm b})
$$
where $\overline{\cJ^m_V}$ is the integral closure of $\cJ_V^m$ in $\cO_X$.
In particular, $\bddK^{[m]}_V$ is always a coherent sheaf.
\endclaim

\noindent
A typical example of what may happen is the Euler vector field linear
space $V\subset T_{\bC^n}\:$:
then the sheaf of holomorphic sections $\cO(V)$ is trivial and generated by
$\varepsilon$, i.e.\ $\cO(V)=\cO_{\bC^n}\varepsilon$, hence its sheaf theoretic
dual is $\cO(V)^*=\cO_{\bC^n}\varepsilon^*$ where $\varepsilon^*$ is the
(unbounded) $1$-form such that $\varepsilon^*\cdot\varepsilon=1$. With the
notation of Proposition~2.2, we have
$$
\eqalign{
\cO(V^*)&=\cO(\Lambda^1V^*)=\bddK_V=\gm_0\cO_{\bC^n}\varepsilon^*=\bddK_V,\cr
\cL_V&=\cO(V^*)^{**}=\cO_{\bC^n}\varepsilon^*=\cO(V)^*,\cr
\cJ_V&=\gm_0,\qquad \overline{\cJ_V^m}=\cJ_V^m=\gm_0^m.\cr}
$$
It is equally important to understand the effect of modifications on the
sheaves $\bddK^{[m]}_V$.

\claim 2.3. Proposition|For any modification 
$\mu:(\widehat X,\widehat V)\to (X,V)$, there are always well defined
injective natural morphisms of rank $1$ sheaves
$$
\bddK_V^{[m]}\hookrightarrow \mu_*\big(\bddK_{\widehat V}^{[m]}\big)
\hookrightarrow \cL_V^{\otimes m}
\leqno(2.3\,{\rm a})
$$
and  the direct image $\mu_*\big(\bddK_{\widehat V}^{[m]}\big)$ 
may only increase when $\mu$ is replaced by a ``higher'' modification
$\widetilde\mu=\mu'\circ\mu:\widetilde X\to\widehat X\to X$
and $\swh{V}=\mu^*V$ by $\swt{V}=\swt{\mu}^*V$, i.e.\ there are injections
$$
\mu_*\big(\bddK_{\widehat V}^{[m]}\big)\hookrightarrow
\widetilde \mu_*\big(\bddK_{\widetilde V}^{[m]}\big)\hookrightarrow
\cL_V^{\otimes m}.
\leqno(2.3\,{\rm b})
$$
We~refer to this property as the {\rm monotonicity principle}.
\endclaim

\plainproof. (a) The existence of the first arrow is seen as follows:
the differential $\mu_*=d\mu:\swh{V}\to \mu^*V$ is smooth, hence bounded
with respect to ambient hermitian metrics on $X$ and $\swh{X}$, and going
to the duals reverses the arrows while preserving boundedness with
respect to the metrics. We thus get an arrow
$$
\mu^*({}^bV^\star)\hookrightarrow{}^b\widehat V^\star.
$$
By taking the top exterior power, followed by the $m$-th tensor product and 
the integral closure of the ideals involved, we get an injective arrow
$\mu^*\big(\bddK_V^{[m]}\big)\hookrightarrow\bddK_{\widehat V}^{[m]}$.
Finally we apply the direct image fonctor $\mu_*$ and the canonical
morphism $\cF\to\mu_*\mu^*\cF$ to get the first inclusion morphism.
The second arrow comes from the fact that $\mu^*\big(\bddK_V^{[m]}\big)$
coincides with $\cL_V^{\otimes m}$ (and with $\det(V^*)^{\otimes m}$) on 
the complement of
the codimension~$2$ set $S=\Sing(V)\cup\mu({\rm Exc}(\mu))$, and the fact
that for every open set $U\subset X$, sections of $\cL_V$ defined on 
$U\ssm S$ automatically extend to $U$ by the Riemann extension theorem,
even without any boundedness assumption.
\smallskip

\noindent (b) Given $\mu':\swt{X}\to\swh{X}$, we argue as in (a) that there is a
bounded morphism $d\mu':\swt{V}\to\swh{V}$.\qed
\medskip

By the monotonicity principle and the strong Noetherian 
property of coherent sheaves, we infer that there exists a maximal 
direct image when $\mu:\swh{X}\to X$ runs over all nonsingular 
modifications of~$X$. The following definition is thus legitimate.

\claim 2.4. Definition|We define the pluricanonical sheaf
$K_V^{[m]}$ of $(X,V)$ to be the inductive limit
$$
K^{[m]}_V:=\lim_{{\scriptstyle\lra\atop\scriptstyle\mu}}
\mu_*\big(\bddK_{\widehat V}^{[m]}\big)=
\max_\mu\mu_*\big(\bddK_{\widehat V}^{[m]}\big)
$$
taken over the family of all modifications $\mu:(\swh{X},\swh{V})\to (X,V)$,
with the trivial $($filtering$)$ partial order. The canonical sheaf $K_V$ 
itself is defined to be the same as $K_V^{[1]}$. By construction, 
we have for every $m\ge 0$ inclusions
$$
\bddK_V^{[m]}\hookrightarrow K_V^{[m]}\hookrightarrow \cL_V^{\otimes m},
$$
and $K_V^{[m]}=\cJ_V^{[m]}\cdot\cL_V^{\otimes m}$ for a certain sequence 
of integrally closed ideals $\cJ_V^{[m]}\subset\cO_X$.
\endclaim

It is clear from this construction that $K^{[m]}_V$ is birationally invariant,
i.e.\ that $K^{[m]}_V=\mu_*(K^{[m]}_{V'})$ for every modification
$\mu:(X',V')\to (X,V)$. One of the most central conjectures
in the theory is~the

\claim 2.5. Generalized Green-Griffiths-Lang conjecture|
Let $(X,V)$ be a projective directed variety. Assume that $(X,V)$ is of
``general type'' in the sense that there exists $m\geq 1$ such that
the invertible sheaf $\mu_m^*(K_V^{[m]})$ is a big line bundle when one takes
a log-resolution $\mu_m$ of the ideal~$\cJ_V^{[m]}$.
Then there should exist a proper algebraic subvariety $Y\subsetneq X$
containing the images $f(\bC)$ of all entire curves $f:\bC\to X$
tangent to~$V$.
\endclaim

The main reason for incorporating the ideals $\cJ_V^{[m]}$ in the
definition of $K_V^{[m]}$ is that the above conjecture would otherwise
be trivially false:
for instance, it is easy to see that a pencil of conics passing
through 4 points in general position in $\bP^2_\bC$ are tangent to a
rank $1$ subspace $V\subset T_{\bP^2}$ such that $\cL_V=\cO(V)^*\simeq
\cO_{\bP^2}(1)$ is ample; however, all leaves are rational curves. Here,
we are in fact more interested here in the dual situation where
a positivity assumption is made for the
line bundle $ \big(\mu_m^*(K_V^{[m]})\big)^*$.

\section{Jet bundles and jet differentials}

\plainsubsection 3.A. Nonsingular case|

Following Green-Griffiths [GrGr79], we consider the bundle $J_kX\to X$
of $k$-jets of germs of para\-metrized curves in~$X$, i.e., the set
of equivalence
classes of holomorphic maps $f:(\bC,0)\to(X,x)$, with the equivalence
relation $f\sim g$ if and only if all derivatives $f^{(j)}(0)=g^{(j)}(0)$
coincide for $0\le j\le k$, when computed in some local coordinate system
of $X$ near~$x$. The projection map $J_kX\to X$ is simply $f\mapsto f(0)$.
If $(z_1\ld z_n)$ are local holomorphic coordinates on an open set
$\Omega\subset X$, the elements $f$ of any fiber $J_kX_x$,
$x\in\Omega$, can be seen as $\bC^n$-valued maps
$$
f=(f_1\ld f_n):(\bC,0)\to\Omega\subset\bC^n,
$$
and they are completetely determined by their Taylor expansion of
order $k$ at~$t=0$
$$
f(t)=x+t\,f'(0)+{t^2\over 2!}f''(0)+\cdots+{t^k\over k!}f^{(k)}(0)+
O(t^{k+1}).
$$
In these coordinates, the fiber $J_kX_x$ can thus be identified with the
set of $k$-tuples of vectors 
$(\xi_1,\ldots,\xi_k)=(f'(0)\ld f^{(k)}(0))\in(\bC^n)^k$.
It follows that $J_kX$ is a holomorphic fiber bundle with typical fiber
$(\bC^n)^k$ over $X$. However, $J_kX$ is not a vector bundle for $k\ge 2$,
because of the nonlinearity of coordinate changes: a coordinate change
$z\mapsto w=\Psi(z)$ on $X$ induces a polynomial
transition automorphism on the fibers of $J_kX$, given by a formula
$$
(\Psi\circ f)^{(j)}=\Psi'(f)\cdot f^{(j)}+\sum_{s=2}^{s=j}{~}
\sum_{j_1+j_2+\cdots+j_s=j}c_{j_1\ldots j_s}\Psi^{(s)}(f)\cdot
(f^{(j_1)}\ld f^{(j_s)})\leqno(3.1)
$$
with suitable integer constants $c_{j_1\ldots j_s}$ (this is easily
checked by induction on~$s$). According to the above philosophy,
we introduce the concept of jet bundle in the
general situation of complex directed manifolds, assuming $V$
nonsingular to start with.

\claim 3.2.~Definition|Let $(X,V)$ be a complex directed manifold.
We define $J_kV\to X$ to be the bundle of $k$-jets of curves
\hbox{$f:(\bC,0)\to X$} that are tangent to $V$, i.e., such that
$f'(t)\in V_{f(t)}$ for all $t$ in a neighborhood of~$0$, together with
the projection map $f\mapsto f(0)$ onto~$X$.
\endclaim

For any point $x_0\in X$, there are local coordinates $(z_1\ld z_n)$ on a
neighborhood $\Omega$ of $x_0$ such that the fibers $(V_z)_{z\in\Omega}$
can be defined by linear equations
$$
V_z=\Big\{v=\sum_{1\le j\le n}v_j{\partial\over\partial z_j}\,;\,
v_j= \sum_{1\le k\le r}a_{jk}(z)v_k~\hbox{\rm for $j=r+1\ld n$}\Big\},
\leqno(3.3)
$$
where $(a_{jk})$ is a holomorphic $(n-r)\times r$ matrix.
Let $f:D(0,R)\to X$ be a curve tangent to $V$ such that $f(D(0,R))\in\Omega$,
and let $(f_1\ld f_n)$ be the components of $f$ in the coordinates
The curve $f$ is uniquely determined by its initial value
$x=f(0)$ and by the first $r$ components $(f_1\ld f_r)$. Indeed, as $f'(t)\in
V_{f(t)}\,$, we can recover the other components by integrating the system
of ordinary differential equations
$$
f_j'(t)=\sum_{1\le k\le r}a_{jk}(f(t))f_k'(t),\qquad r+1\leq j\leq n,\leqno(3.4)
$$
on a neighborhood of~$0$, with initial data $f(0)=x$. As a consequence,
$J_kV$ is actually a subbundle of~$J_kX$. In fact, by using (3.4), we
see that the fibers $J_kV_x$ are parametrized by 
$$
\big((f_1'(0)\ld f_r'(0));(f_1''(0)\ld f_r''(0));\ldots;
(f_1^{(k)}(0)\ld f_r^{(k)}(0))\big)\in(\bC^r)^k
$$
for all $x\in\Omega$, hence $J_kV$ is a locally trivial
$(\bC^r)^k$-subbundle of~$J_kX$. Alternatively, we can pick a
local holomorphic connection $\nabla$ on $V$ such that for
any germs $w=\sum_{1\leq j\leq n}w_j{\partial\over \partial z_j}
\in\cO(T_{X,x})$ and $v=\sum_{1\leq\ell\leq r}v_\ell e_\ell\in\cO(V)_x$
in a local trivializing frame $(e_1,\ldots,e_r)$ of~$V_{\restriction \Omega}$
we have
$$
\nabla_wv(x)=\sum_{1\leq j\leq n,\,1\leq\ell\leq r}
w_j{\partial v_\ell\over\partial z_j}e_\ell(x)+
\sum_{1\leq j\leq n,\,1\leq\ell,\mu\leq r}
\Gamma_{j\ell}^\mu(x)w_jv_\ell\,e_\mu(x).
\leqno(3.5)
$$
We can of course take the (unique) frame $(e_\ell)_{1\leq\ell\leq r}$ in
$V$ such that $\partial/\partial z_\ell$ is the projection of $e_\ell$ on
the first $r$ coordinates, and the trivial connection $\nabla^0$ given by the 
zero Christoffel symbolds $\Gamma=0$ with respect to this frame, but any
other holomorphic connection $\nabla$ is acceptable. One 
then obtains a trivialization $J^kV_{\restriction \Omega}\simeq
V^{\oplus k}_{\restriction \Omega}$ by considering
$$
J_kV_x\ni f\mapsto(\xi_1,\xi_2,\ldots,\xi_k)=(\nabla f(0),\nabla^2f(0),\ldots,\nabla^k f(0))\in V_x^{\oplus k}
\leqno(3.6)
$$
and computing inductively the successive derivatives $\nabla f(t)=f'(t)$ and
$\nabla^sf(t)$ via
$$
\nabla^sf=(f^*\nabla)_{d/dt}(\nabla^{s-1}f)
=\sum_{1\leq\ell\leq r}
{d\over dt}\Big(\nabla^{s-1}f\Big)_\ell e_\ell(f)+
\sum_{1\leq j\leq n,\,1\leq\ell,\mu\leq r}
\Gamma_{j\ell}^\mu(f)f'_j\Big(\nabla^{s-1}f\Big)_\ell e_\mu(f).
$$
This identification depends of course on the choice of $\nabla$ and cannot 
be defined globally in general (unless we are in the rare situation
where $V$ has a global holomorphic connection). Now,
we consider the natural $\bC^*$-action on $J^kV$ that maps a $k$-jet
$t\mapsto f(t)$ to $\lambda\cdot f(t):=f(\lambda t)$, $\lambda\in\bC^*$.
Since $\nabla^s(\lambda\cdot f)(t)=\lambda^s\nabla^sf(t)$, the $\bC^*$
action is described
in coordinates by
$$
\lambda\cdot(\xi_1,\xi_2,\ldots,\xi_k)=
(\lambda \xi_1,\lambda^2\xi_2,\ldots,\lambda^k\xi_k),\quad \xi_s=\nabla^sf(0).
\leqno(3.7)
$$
Following [GrGr79], we introduce the bundle $E^\GG_{k,m}V^*\to X$ 
of polynomials $P(x\,;\xi_1,\ldots,\xi_k)$ that are homogeneous
on the fibers of $J_kV$ of weighted degree $m$ with respect to the
$\bC^*$ action, i.e.
$$
P(x\,;\lambda\xi_1,\ldots,\lambda^k\xi_k)=\lambda^m P(x\,;\xi_1,\ldots\xi_k),
\leqno(3.8)
$$
in other words they are polynomials of the form
$$
P(x\,;\xi_1,\ldots\xi_k)=
\sum_{|\alpha_1|+2|\alpha_2|+\cdots+k|\alpha_k|=m}
a_{\alpha_1\ldots\alpha_k}(x)\,\xi_1^{\alpha_1}\xi_2^{\alpha_2}\cdots
\xi_k^{\alpha_k}\leqno(3.9)
$$
where $\xi_s=(\xi_{s,1},\ldots,\xi_{s,r})\in\bC^r\simeq V_x$ and
$\xi_s^{\alpha_s}=\xi_{s,1}^{\alpha_{s,1}}\ldots\xi_{s,r}^{\alpha_{s,r}}$,
$|\alpha_s|=\sum_{1\leq j\leq r}\alpha_{s,j}$. Sections of the
sheaf $\cO(E^\GG_{k,m}V^*)$ can also be viewed as algebraic differential
operators acting on germs of curves $f:(\bC,0)\to X$ tangent to~$V$, by putting
$$
P(f)(t)=\sum_{|\alpha_1|+2|\alpha_2|+\cdots+k|\alpha_k|=m}
a_{\alpha_1\ldots\alpha_k}(f(t))\,(\nabla f(t))^{\alpha_1}(\nabla^2 f(t))^{\alpha_2}\cdots
(\nabla^kf(t))^{\alpha_k}\leqno(3.9')
$$
where the $a_{\alpha_1\ldots\alpha_k}(x)$ are holomorphic in~$x$.
With the graded algebra bundle $E^\GG_{k,\bu}V^*=\bigoplus_mE^\GG_{k,m}V^*$
we associate an analytic fiber bundle
$$X_k^\GG:=\Proj(E^\GG_{k,\bu}V^*)=(J_kV\ssm\{0\})/\bC^*\leqno(3.10)$$
over $X$, which has weighted projective spaces
$\bP(1^{[r]},2^{[r]}\ld k^{[r]})$ as fibers; here $J_kV\ssm\{0\}$ is the set
of nonconstant jets of order~$k$.
As such, it possesses a tautological sheaf $\cO_{X^\GG_k}(1)$ [the reader
should observe however that $\cO_{X^\GG_k}(m)$ is invertible only when
$m$ is a multiple of $\lcm(1,2,\ldots,k)$].

\claim 3.11. Proposition|By construction, if $\pi_k:X_k^\GG\to X$ is the
natural projection, we have the direct image formula
$$
(\pi_k)_*\cO_{X^\GG_k}(m)=\cO(E_{k,m}^\GG V^*)
$$
for all $k$ and $m$.
\endclaim

\plainsubsection 3.B. Singular case|

When $V$ has singularities and $X$ is nonsingular, we simply consider
the inclusion morphism $(X,V)\to (X,T_X)$ into the absolute directed structure.
This yields a morphism $J_kV_{|X'}\to J_kT_{X|X'}$ in restriction to the Zariski
open set $X'=X\ssm\Sing(V)$, and we define $J_kV=\overline{J_kV_{|X'}}$
to be the closure of $J_kV_{|X'}$ in $J_kT_X$. It is then easy to see
that $J_kV$ is an analytic subset of $J_kT_X$, hence we get
an inclusion morphism $J_kV\hookrightarrow J_kT_X$ over~$X$, which also induces
an inclusion
$$
X_k^\GG\hookrightarrow X_k^{\abs,\GG}
\leqno(3.12)
$$
of $X_k^\GG=(J_kV\ssm\{0\})/\bC^*$ into the absolute Green-Griffiths bundle
$X_k^{\abs,\GG}=(J_kT_X\ssm\{0\})/\bC^*$. In analogy with our concept of
canonical sheaf, it is natural to introduce the following definitions.

\claim 3.13. Theorem and definition|Let $p_k:J^kV\ssm\{0\}\to X^\GG_k$
be the natural projection. The sheaf ${}^b\cO_{X^\GG_k}(m)$ $($here, the $b$
means ``locally bounded'' sections$\,)$ is defined as follows: for any open set
$U\subset X_k^\GG$, the space of sections ${}^b\cO_{X^\GG_k}(m)(U)$
consists of holomorphic operators $F(x\,;\xi_1,\ldots,\xi_k)$ on
the conical open set
$$
p_k^{-1}(U)\cap J^kV_{|X'}\subset J^kV_{|X'}
\subset J^kT_X
$$
that are homogeneous of degree $m$ with respect to the $\bC^*$-action, and
are locally bounded with
respect to a smooth hermitian metric $h_X$ on $T_X$. Namely,
if $\nabla$ is a local holomorphic connection
on $T_{X|\Omega}$ and $(\xi_1,\ldots,\xi_k)$ are the components of a $k$-jet
computed with respect to $\nabla$, we require that
$$
|F(x\,;\xi_1,\ldots,\xi_k)|\leq C(U_c)\bigg(\sum_{1\leq s\leq k}\Vert\xi_s\Vert_{h_X}^{1/s}\bigg)^m\leqno(3.13^*)
$$
on $p_k^{-1}(U_c)$, for every relatively compact open subset $U_c\compact U
\cap\pi_k^{-1}(\Omega)$. Then ${}^b\cO_{X^\GG_k}(m)$ is a rank~$1$
coherent analytic sheaf, and is independent of the choice of $h_X$ and $\nabla$.
\endclaim

\claim 3.14. Definition|With the above notation, the sheaf 
${}^b\cO(E_{k,m}^\GG V^*)$ is the analytic sheaf on $X$ whose spaces of
sections are polynomial differential operators $P(x\,;\xi_1,\ldots,\xi_k)$ in
$$
{}^b\cO(E_{k,m}^\GG V^*)(U)\subset\cO(E_{k,m}^\GG V^*)(U\cap X'),
$$
i.e. polynomial functions $P=F$ that satisfy inequality $(3.13^*)$
on $\pi_k^{-1}(U_c\cap X')$, for all open sets $U_c\compact U\subset X$.
In other words, we have
$$
(\pi_k)_*{}^b\cO_{X^\GG_k}(m)={}^b\cO(E_{k,m}^\GG V^*).
\leqno(3.14^*)
$$
\endclaim

\proof That homogeneous functions on $J^kV_{|U}$ must be polynomials on
the fibers is a trivial fact (using e.g.\ power series expansions in terms
of $(\xi_1,\ldots,\xi_k)$). The coherence of ${}^b\cO_{X^\GG_k}(m)$
is a simple consequence
of the fact that we have a restriction morphism
$$
\cO_{X^{\abs,\GG}_k}(m)_{|X^\GG_k}\longrightarrow{}^b\cO_{X^\GG_k}(m)
\leqno(3.15)
$$
and that, almost by definition, $\bigoplus_m{}^b\cO_{X^\GG_k}(m)$ consists of
taking the normalization of the image of the graded ring of
sections; this shows again that our concepts are purely algebraic, in
spite of the analytic definition that was given in 3.13.
The coherence of ${}^b\cO(E_{k,m}^\GG V^*)$ then follows by the
direct image theorem. It is very important to observe that in condition
$(3.13^*)$ one must refer to a metric $h_X$ and a connection $\nabla$
on the ambient bundle $T_X$, and not to a holomorphic connection on $V$,
since any such connection would be unbounded near $\Sing(V)$, leading to the
failure of (3.15).\qed

It is also clear form the definitions that our sheaf of jet differentials
coincides with the already defined vector bundle on $X'=X\ssm\Sing(V)\,$:
$$
\Big({}^b\cO(E_{k,m}^\GG V^*)\Big)_{|X\ssm\Sing(V)}=
\cO(E_{k,m}^\GG V^*_{|X\ssm\Sing(V)}).\leqno(3.16)
$$
When $X$ itself has singularities, one can locally embed $X$ in a smooth
ambient variety $Z$ and refer to smooth hermitian metrics on $Z$, taking
bounded sections on a Zariski open set where both $X$ and $V$ are smooth.
We will not consider this case much further and leave details to the reader.

\section{Morse inequalities, in the smooth and singular
contexts}

\plainsubsection 4.A. Smooth case|

The main purpose of holomorphic Morse inequalities is to provide
estimates of cohomology groups with values in high tensor powers
of a given line bundle $L$, once a hermitian metric $h$ on $L$ is given.
We denote by $\Theta_ {L,h}=-{i\over 2\pi}\ddbar\log h$ the $(1,1)$-curvature
form of $h$, and for any $(1,1)$-form $u(z)=i\sum_{1\leq j,k\leq n}u_{jk}(z)\,
dz_j\wedge d\overline z_k$, we define its $q$-index set $X(u,q)$ to be
the open set
$$
X(u,q)=\big\{z\in X\,;\;u(z)~\hbox{has signature }~(n-q,q)\big\}
\leqno(4.1)
$$
(so that $q$ is the number of negative eigenvalues of $u(z)$). The 
following statement summarizes the main results of [Dem85].

\claim 4.2. Holomorphic Morse inequalities for smooth metrics|Let $X$
be a compact complex $n$-dimensional manifold, $E\to X$ a holomorphic 
vector bundle of
rank $r$, and $L$ a holomorphic line bundle equipped with a smooth hermitian
metric $h$. The dimensions $h^q(X,E\otimes L^m)$ of cohomology groups of
the tensor powers $E\otimes L^m$ satisfy the following asymptotic estimates as 
$m\to +\infty\,:$\vskip4pt
\noindent {\rm(4.2~WM)} Weak Morse inequalities$\,:$
$$h^q(X,E\otimes L^m)\le r {m^n\over n!}\int_{X(\Theta_{L,h},q)} (-1)^q \Theta_{L,h}^n + o(m^n)~.$$
\noindent {\rm(4.2~SM)} Strong Morse inequalities$\,:$
$$\sum_{0\le j\le q} (-1)^{q-j}h^j(X,E\otimes L^m) \le r {m^n\over n!}
\sum_{0\leq j\leq q}\int_{X(\Theta_{L,h},j)}(-1)^{q-j}\Theta_{L,h}^n+o(m^n)~.$$
\noindent {\rm(4.2~RR)} Asymptotic Riemann-Roch formula$\,:$
$$\chi(X,E\otimes L^m) := \sum_{0\le j\le n} (-1)^jh^j(X,E\otimes L^m)
= r{m^n\over n!}\int_X \Theta_{L,h}^n + o(m^n)~.$$
\endclaim

In fact, the strong Morse inequality implies the weak form (by adding 
the inequalities for $q$ and $q-1$), and the asymptotic Riemann-Roch formula
(by taking $q=n$ and $q=n+1$). By adding the strong Morse inequalities for
$q+1$ and $q-2$, one also gets the lower bound
$$\leqalignno{
h^q(X,E\otimes L^m)&\geq h^q(X,E\otimes L^m)-
h^{q-1}(X,E\otimes L^m)-h^{q+1}(X,E\otimes L^m)\cr
&\geq
 r {m^n\over n!}\sum_{q-1\leq j\leq q+1}
\int_{X(\Theta_{L,h},j)} (-1)^{q-j}\Theta_{L,h}^n - o(m^n),&(4.3_q)\cr}
$$
and especially, for the important case $q=0$, the lower bound
$$
h^0(X,E\otimes L^m)\geq 
\int_{X(\Theta_{L,h},0)} \Theta_{L,h}^n -
\int_{X(\Theta_{L,h},1)} \Theta_{L,h}^n - o(m^n).
\leqno(4.3_0)
$$

\plainsubsection 4.B. Case of metrics with $\bQ$-analytic singularities|

The above estimates are volume estimates, and therefore are not sensitive to
modifications. In parti\-cular, we have the following easy lemma.

\claim 4.4. Lemma|Let $\cE$ be a coherent analytic sheaf
and $(L,h)$ a hermitian line bundle on a reduced 
irreducible compact complex space~$Z$. Then
\plainitem{\rm(i)} If $Y=\Supp(\cE)$ has at most $p$-dimensional
irreducible components, we have
$$
h^q(Z,\cE\otimes \cO(L^m))=O(m^p).
$$
\plainitem{\rm(ii)} If $n=\dim Z$ and $r$ is the generic rank of $\cE$, 
the Morse inequalities are still valid with $X$ replaced by $Z$ and 
the locally free sheaf $\cO(E)$ replaced by $\cE$.
\vskip0pt
\endclaim

\proof We prove $4.4\,{\rm(i)}_{p\leq N}$ and $4.4\,{\rm(ii)}_{n\leq N}$ 
by induction on $N$, everything being obvious for $N=0$.

If $Y$ is not irreducible and $Y=\bigcup Y_j$ is 
the decomposition in irreducible components, let $\cE_j$ be the 
sheaf of sections of $\cE$ that
vanish on all components $Y_k$, $k\neq j$. Then we have an exact sequence
$$
0\to\bigoplus\cE_j\to \cE\to\cF\to 0
$$
where $\cF$ is supported on $Y'=\bigcup_{j\neq k}Y_j\cap Y_k$, $\dim Y\leq p-1$.
If we use the corresponding exact sequence and argue by induction on $p$, 
we see that it is sufficient to check $4.4\,{\rm(i)}_p$ when $Y$ is irreducible.
If $\cI_Y$ is the reduced ideal sheaf of $Y$, we 
have $\cI_Y^k\cE=0$ for some $k\in\bN^*$, and we thus get a 
decreasing filtration $\cE_\ell=\cI_Y^\ell\cE$ of $\cE$ such that
$\cE_\ell/\cE_{\ell+1}$ can be viewed as a coherent sheaf on the reduced space
$Y_{\rm red}$, whose structure sheaf is $\cO_X/\cI_Y$. Then, by taking
$Z=Y$ and exploiting the exact sequences
$0\to\cE_{\ell+1}\to\cE_\ell\to\cE_\ell/\cE_{\ell+1}\to 0$, we see 
that $4.4\,{\rm(ii)}_{n\leq N}$ implies $4.4\,{\rm(i)}_{p\leq N}$. The last
part of the proof consists in showing that Theorem~$4.2_{\,n\leq N}$ and
$4.4\,{\rm(i)}_{p\leq N-1}$ imply $4.4\,{\rm(ii)}_{n\leq N}$, and
for this, it is enough to consider the case where $\dim Z=n=N$.

In that case, the Hironaka desingularization theorem implies that there exists
a modification $\mu:X\to Z$ such that $X$ is nonsingular and
$\cF=\mu^*\cE$ is locally free modulo torsion, so that we have an exact 
sequence
$$
0\to \cF_\tors\to\cF\to \cF/\cF_\tors\to 0
$$
where $\cF/\cF_\tors$ is locally free (i.e.\ a vector bundle on $X$). Therefore,
holomorphic inequalities 4.2 can be applied to the groups $H^q(X,
\cF/\cF_\tors\otimes\cO(\mu^*L^m))$, and since $Y=\Supp(\cF_\tors)$ has
dimension $p\leq N-1$, part $4.4\,{\rm(i)}_{p\leq N-1}$ of the Lemma 
shows that the groups $H^q(X,\cF\otimes\cO(\mu^*L^m))$ also satisfy
holomorphic Morse inequalities. Finally, we use the Leray spectral
sequence. It yields a convergent spectral sequence
$$
H^p(Z,R^q\mu_*(\cF)\otimes\cO(L^m))\Rightarrow 
H^{p+q}(X,\cF\otimes\cO(\mu^*L^m))
$$
and we know that $R^q\mu_*(\cF)$ is supported for $q\geq 1$ on an analytic
subset $Y'\subsetneq Z$, and that the morphism 
$\cE\to\mu_*\mu^*\cE=R^0\mu_*(\cF)$ is an isomorphism outside of 
codimension~$1$. From this we conclude that holomorphic Morse inequalities
for the $H^q(X,\cF\otimes\cO(\mu^*L^m))$ (which we already know), imply
holomorphic Morse inequalities for $H^q(Z,\cE\otimes\cO(L^m))$,
thanks to $4.4\,{\rm(i)}_{n\leq N-1}$ applied on~$Z$.
\qed

In his PhD thesis, Bonavero [Bon93] extended Morse inequalities to the
case of singular hermitian metrics. We state here a variant that
allows more general singularities (but in fact, everything can be reduced
to the smooth case 4.2 by means of a desingularization, thus all versions are 
in fact equivalent).

\claim 4.5. Definition|Let $Z$ be an irreducible and reduced complex space 
and $\cL$ a torsion free rank~$1$ sheaf on $Z$. A hermitian metric $h$
on $\cL$ with $\bQ$-analytic singularities is a hermitian metric $h$ defined
on a dense open set $Z'\subset Z_{\reg}$ where $\cL$ is invertible, with the following property$\,:$ there exists a smooth modification $\mu:X\to Z$
such that $\mu^*\cL$ is invertible and the pull-back metric $\mu^*h$ has
normal crossing $\bQ$-divisorial singularities, i.e.\ if $g$ is a
local generator of $\mu^*\cL$ on a small open set $U\subset X$, we have
$$
\log|g|^2_{\mu^*h}=\psi-\sum\lambda_j\log|z_j|^2
$$
for some holomorphic coordinates $(z_j)$ on $U$, $\lambda_j\in\bQ\ssm\{0\}$
and $\psi\in C^\infty(U)$. We will say that $h$ has
$\bQ_+$-analytic singularities if one can take all
$\lambda_j\in\bQ_+\ssm\{0\}$ $($or an empty sum$\,)$.
\endclaim

(Of course the normal crossing hypothesis is not necessary, since it can always
be achieved a posteriori by an application of the Hironaka desingularization
theorem). We define the singular set $S=\Sing(\mu^*h)$ to be the common
zero set $\{z_j=0\}$ for the
coordinates $z_j$ involved above, and $\Sing(h)$ to be the union of $\mu(S)$
and of the closed analytic subset of $Z$ where $\cL$ is not invertible.

\claim 4.6. Definition|Let $Z$ be an irreducible and reduced complex space,
and $\cL$ a torsion free rank~$1$ sheaf on $Z$ equipped with a
hermitian metric $h$ with $\bQ$-analytic singularities. We define
the sheaf of $h$-bounded sections ${}^b\cO(\cL^m)_h$ to be the sheaf
of germs of meromorphic sections $\sigma$ of $\cL^m$ such that
$|\sigma|_h$ is bounded $($just consider the restriction to
$U\cap Z'$ where $U$ is a neighborhood of a given point $z_0\in Z$ and
$Z'=Z\ssm\Sing(h)\,)$.
\endclaim

It follows from the definitions that ${}^b\cO(\cL^m)_h$ is a coherent sheaf
of rank~$1\,$; in fact, with the notation of Definition~4.5, it is the
direct image by $\mu$ of $(\mu^*\cL)^m\otimes\cO_X
(-\lceil mD\rceil)$
where $D=\sum \lambda_jD_j$ is the $\bQ$-divisor of $X$ given by
$D_j=\{z_j=0\}$, and $\lceil mD\rceil=\sum \lceil m\lambda_j\rceil D_j$
is the round up. Especially, if $Z$ is normal, $\cL$ invertible and $h$
has $\bQ_+$-analytic singularities (i.e.\ $D\geq 0$), then we have
\hbox{${}^b\cO(\cL^m)_h\subset \cO(\cL^m)$}. Otherwise, we may get some
meromorphic sections in ${}^b\cO(\cL^m)_h$ that are actually not holomorphic.

\claim 4.7. Holomorphic Morse inequalities for singular metrics|Let $Z$ be
an irreducible and reduced complex space, $\cL$ a torsion free rank~$1$
sheaf on $Z$ equipped with a hermitian metric $h$ with $\bQ$-analytic
singularities, and $\cE$ be a coherent sheaf of generic rank~$r$.
Then, for $Z'=Z\ssm\Sing(h)$ and all $q=0,1,\ldots,n=\dim Z$, we have
estimates
$$
\eqalign{
\sum_{j=q-1,q,q+1}
r{m^n\over n!}\int_{Z'(\Theta_{L,h},j)}
&(-1)^{q-j}\Theta_{\cL,h}^n-o(m^n)\cr
&\leq h^q(Z,\cE\otimes{}^b\cO(\cL^m)_h)\leq r{m^n\over n!}
\int_{Z'(\Theta_{L,h},q)}(-1)^q\Theta_{\cL,h}^n+o(m^n).\cr}
$$
\endclaim

\proof We use a smooth modification $\mu:X\to Z$ such that $\mu^*\cL$
satisfies the requirements of Definition 4.5 and $\mu^*\cE$ is locally
free modulo torsion. The proof of Lemma~4.4 then shows that the torsion
of $\cF=\mu^*\cE$ can be neglected. Again, the conclusion follows from
the Leray spectral sequence and the smooth case of Morse inequalities
applied to the sequence of invertible line bundles
$$
p\mapsto\cG\otimes\widetilde\cL^p\quad\hbox{on $X$},
$$
where
$$
\cG=\cF/\cF_\tors\otimes\mu^*\cL^r
\otimes\cO_X(-\lceil rD\rceil),\quad
\widetilde\cL=\mu_\star\cL^a\otimes\cO_X(-\lceil aD\rceil),
$$
$a\in\bN^*$ is a denominator of the $\bQ$-divisor $D$ and $m=ap+r$,
$0\leq r<a$. It follows that there is a morphism
$$
\cE\otimes{}^b\cO(\cL^m)_h) \to \mu_*(\cG\otimes\widetilde\cL^p)
$$
whose kernel and cokernel are supported on subvarieties
$Y,Y'\subsetneq Z$, and by definition $\mu^*h^a$ induces a smooth metric
on $\widetilde\cL$. As a consequence
$\Theta_{\widetilde\cL,\mu^*h^a}=a\,\mu^*\Theta_{\cL,h}$ on
$\mu^{-1}(Z\ssm\Sing(h))$ and
$$
\int_{X(\Theta_{\widetilde\cL,\mu^*h^a},q)}(-1)^q\Theta_{\widetilde\cL,\mu^*h^a}
^{\raise3pt\hbox{$\scriptstyle n$}} =
a^n\int_{Z'(\Theta_{L,h},q)}(-1)^q\Theta_{\cL,h}^n.\eqno\square
$$

\claim 4.8. Remark|{\rm In Bonavero's thesis [Bon93], $Z$ is a manifold,
$\cE=\cO_Z(E)$ and $\cL=\cO_Z(L)$ are assumed to be locally free,
$h=e^{-\varphi}$ is a singular hermitian metric with
$\bQ_+$-analytic singularities and the cohomology groups involved are
the groups
$$
H^q(Z,\cO_Z(E\otimes L^m)\otimes\cI(h^m))
$$
where $\cI(h^m)$ are the $L^2$ multiplier ideal sheaves 
$$
\cI(h^m)=\cI(k\varphi)=\big\{f\in\cO_{Z,x},\;\;\exists V\ni x,~
\int_V|f(z)|^2e^{-m\varphi(z)}d\lambda(z)<+\infty\big\}.
$$
Here, we have in fact replaced the $L^2$ condition by a $L^\infty$ condition.
When pulling-back to $X$ via a modification $\mu:X\to Z$ resolving the
singularities of $h$ into divisorial singularities, the difference is just
a twist by the relative canonical bundle $K_{X/Z}$ and the use of the
round down $\lfloor m D\rfloor$ instead of~$\lceil mD\rceil$. These differences
are ``bounded'' and thus do not make any change in the estimates produced
by the Morse inequalities. For the same reason, we could even allow the
singularity $D$ of $\mu^*h$ to be a $\bR$-divisor, although periodicity
would be lost in the round up process; one then needs the fact
that Morse inequalities are ``uniform'' when $E$ remains in a bounded
family of vector bundles, which follows easily from the analytic
proof.}
\endclaim

\section{Cohomology estimates for sheaves of jet differentials
and their duals}

On a directed variety $(X,V)$, is may be necessary to allow $V$
to have singular hermitian metrics, even when $(X,V)$ is smooth: indeed,
this leads to more comprehensive curvature conditions, since one can
e.g.\ relax the ampleness conditions and consider instead big line bundles.
However, it is also useful to consider situations where $V$ is singular.
According to the philosophy of \S$\,$4, we have to explain what are
the sheaves involved in the presence of such singularities.

\claim 5.1. Definition|Let $(X,V)$ be a directed variety, $X$ being nonsingular.
\vskip3pt
\plainitem{\rm(a)} A singular hermitian metric on a linear
subspace $V\subset T_X$ is a metric $h$ on the fibers of $V$ such that
the function $\log h:\xi\mapsto\log|\xi|_h^2$ is locally integrable 
on the total space of $V$.
\plainitem{\rm(b)} A singular metric $h$ on $V$ will be said to have
$\bQ$-analytic singularities if $h$ can be written as $h=e^\varphi
(h_X)_{|V}$ where $h_X$ is a smooth positive definite hermitian metric
on $T_X$ and $\varphi$ is a weight with $\bQ$-analytic singularities, i.e.\
one can find a modification $\mu:\widetilde X\to X$ such that locally
$\varphi\circ\mu=\sum\lambda_j\log|g_j(z)|+\psi(z)$, where
$\lambda_j\in\bQ\ssm\{0\}$, $\psi\in C^\infty$ and the $g_j$ are
holomorphic functions; one can then further assume that
$D=\sum\lambda_j D_j$, $D_j=\{g_j=0\}$, is a simple normal crossing divisor
on~$\smash{\widetilde X}$. We define $\Sing(h)$ to be the union
$\mu(\Supp(D))\cup\Sing(V)$.
\vskip0pt
\endclaim

A singular metric $h$ on $V$ can also be viewed as a singular
hermitian metric on the tauto\-logical line bundle $\cO_{P(V)}(-1)$
on the projectivized bundle
$P(V)=V\ssm\{0\}/\bC^*$, and therefore its dual metric $h^{-1}$ defines a 
curvature current $\Theta_{\cO_{P(V)}(1),h^{-1}}$ of type $(1,1)$ on
$P(V)\subset P(T_X)$, such that
$$
p^*\Theta_{\cO_{P(V)}(1),h^{-1}} ={\ii\over 2\pi}\ddbar\log h,\qquad
\hbox{where $p:V\ssm\{0\}\to P(V)$}.\leqno(5.2)
$$

\claim 5.3. Remark|{\rm In [Dem11], [Dem12], [Dem15], we introduced the concept
of an admissible metric $h$ on~$V$, which is closely related to the
concept of a metric with $\bQ_+$-analytic singularities (in the sense
that the divisor of singularities $D$ is nonnegative);
then $\log h$ is quasi-plurisubharmonic $($i.e.\ psh
modulo addition of a smooth function$)$ on the total space of $V$, hence
$$
\Theta_{\cO_{P(V)}(1),h^{-1}}\ge -C\omega\leqno(5.3\,{\rm a})
$$
for some smooth positive $(1,1)$-form $\omega$ on $P(V)$ and some
constant $C>0\;$; if $h$ has $\bQ$-analytic singularities, we can choose
a smooth modification $\widetilde X\to X$ such that $\cO(\mu^*V)$ is
locally free,
the injection $\cO(\mu^*\det V)\hookrightarrow\cO(\mu^*\Lambda^rT_X)$
vanishes along an invertible ideal $\cO(-\Delta)\subset \cO_{\widetilde X}$
and $\mu^*\log(h/h_X)$ has divisorial singularities given by a
$\bQ$-divisor $D$. Then
$$
\Theta_{\cO_{P(\mu^*V)}(1),\mu^*h^{-1}}=[D+\Delta]+\beta\leqno(5.3\,{\rm b})
$$
for some smooth $(1,1)$-form $\beta$ on $P(\mu^*V)$. Hence if $D+\Delta\geq 0$
(and especially if $D\geq 0$), we still have
$$
\Theta_{\cO_{P(\mu^*V)}(1),\mu^*h^{-1}}\geq -C\widetilde\omega\leqno(5.3\,{\rm c})
$$
for any smooth positive $(1,1)$-form $\widetilde\omega$ on~$P(\mu^*V)$.
}
\endclaim

\claim 5.4. Bounded sections|
If $h$ is a singular metric with $\bQ$-analytic singularities
on $(X,V)$ and~$h^*$ the corresponding
dual metric on $V^*$, we consider the Zariski open set
$X'=X\ssm(\Sing(V)\cup\Sing(h))$ and define$\,:$
\plainitem{\rm (a)} ${}^b\cO(V^*)_{h^*}$ to be the sheaf of germs
of holomorphic sections of $V^*_{|X'}$ which are $h^*$-bounded near every point
of~$X\,;$
\plainitem{\rm (b)} the $h$-bounded pluricanonical sequence to be
$$
\eqalign{
\bddK_{V,h}^{[m]}={}&\hbox{sheaf of germs of holomorphic sections of 
$(\det V^*_{|X'})^{\otimes m}=(\Lambda^rV^*_{|X'})^{\otimes m}$}\cr
\noalign{\vskip-2pt}
&\hbox{which are $\det h^*$-bounded},\cr}
$$
so that $\bddK_V^{[m]}:=\bddK_{V,h_X}^{[m]}$ according to Def.~$2.1$.
\plainitem{\rm (c)} for $U\subset X$ open, the space of sections
${}^b\cO(E^\GG_{k, m}V^*_h)(U)$ of ${}^b\cO(E^\GG_{k, m}V^*_h)$
consists of functions $P(x\,;\xi_1,\ldots,\xi_k)$ that are holomorphic
in $x\in U\cap X'$ and polynomial in the $\xi_j$'s, satisfying
upper bounds of the form
$$
\big|P(z\,;\,\xi_1,\ldots,\xi_k)\big|\le 
C(U_c)\bigg(\sum_{1\leq s\leq k}\Vert \xi_s\Vert_h^{1/s}\bigg)^m,\quad
x\in U_c\compact U.
$$
\vskip0pt\noindent
We then have a direct image formula
$$
{}^b\cO(E^\GG_{k, m}V^*_h)=(\pi_k)_*{}^b\cO_{X_k^\GG}(m)_h.\leqno{\rm (d)}
$$
where ${}^b\cO_{X_k^\GG}(m)_h$ denotes the sheaf of $h$-bounded
sections of ${}^b\cO_{X_k^\GG}(m)$.
\endclaim

If the divisor $D$ of singularities of $h$ is${}\geq 0$,
we have $h=e^\varphi h_X\leq Ch_X$ locally, hence
$$
\bddK_{V,h}^{[m]}\subset\bddK_V^{[m]}\quad\hbox{and}\quad
{}^b\cO(E^\GG_{k, m}V^*_h)\subset{}^b\cO(E^\GG_{k, m}V^*).
\leqno(5.5)
$$
On the other hand, if $D\leq 0$, reversed inequalities and
inclusions hold, e.g.\
$\bddK_{V,h}^{[m]}\supset\bddK_V^{[m]}$, therefore
$$
(\bddK_{V,h}^{[m]})^*\subset(\bddK_V^{[m]})^*\quad\hbox{and}\quad
\cO({}^bE^\GG_{k, m}V^*_h)^*\subset\cO({}^bE^\GG_{k, m}V^*)^*
\leqno(5.5^*)
$$
for the dual $\cO_X$-modules.

\claim 5.6. Morse integral estimates|Let $(X,V)$ be a directed
variety, where $X$ is a nonsingular compact complex manifold.
Fix a singular hermitian metric $h_V$ on $V$ and denote by
$$
\Theta_{V,h_V}={i\over 2\pi}\sum_{1\leq i,j\leq n,1\leq\alpha,\beta\leq r}
c_{ij\alpha\beta}(z)\,dz_i\wedge d\overline z_j\otimes e_\alpha^*\otimes e_\beta
$$ the curvature tensor of $V$ with respect to an $h_V$-orthonormal frame
$(e_\alpha)$, in the sense of currents. Let $(F,h_F)$ be a singular
$\bQ$ hermitian line bundle on $X$, and let
$$
\eta(z)=\Theta_{\det(V^*),\det h_V^*}+\Theta_{F,h_F}=-\Tr_{\End(V)}\Theta_{V,h_V}
+\Theta_{F,h_F}.
\leqno(5.6\,{\rm a})
$$
Assume that $h=e^\varphi(h_X)_{|V}$ and $h_F$ are metrics 
with $\bQ$-analytic singularities and let $\Sigma$ be their joint singular
set $\Sigma=\Sing(h_V)\cup\Sing(h_F)$. Finally, equip the tautological bundle
$\cO_{X_k^\GG}(-1)$ with the induced metric $h_{V,k,\varepsilon}$
such that
$$
\Vert\xi\Vert^2_{h_{V,k,\varepsilon}}=e^\varphi
\Bigg(\sum_{1\leq s\leq k}\varepsilon_s|\xi_s|_{h_X}^{2p/s}\Bigg)^{1/p},
\quad p=\lcm(1,2,\ldots,k),
\leqno(5.6\,{\rm b})
$$
where $\xi=f_{[k]}(0)$ is the $k$-jet of an integral germ of curve
$f$ in $(X,V)$, $\xi_s=(\nabla_0)^sf(0)$ with respect to a global smooth
connection $\nabla_0$ on $T_X$ and
$1=\varepsilon_1\gg\varepsilon_2\gg\ldots\gg\varepsilon_k>0$. We
consider on $X_k^\GG$ the rank $1$ sheaf
$$
\cL_k=\cO_{X_k^\GG}(1)\otimes
\pi_k^*\bigg(\cO\bigg({1\over kr}\Big(1+{1\over 2}+\ldots+{1\over k}
\Big)F\bigg)\bigg)
\leqno(5.6\,{\rm c})
$$
$($it is just a rank $1$ ``$\bQ$-sheaf'', and we actually have to take a power
$\cL_k^m$ with $m$ sufficiently divisible to get a genuine rank $1$ sheaf$\,)$,
equipped with the metric $h_{V,F,k,\varepsilon}$ induced by
$(h_{V,k,\varepsilon})^{-1}$ and $h_F$.
Then, for $m\gg\varepsilon_k^{-1}\gg k\gg 1$, the $q$-index Morse integral
of $(\cL_k,h_{V,F,k,\varepsilon})$ on $X_k^\GG$ is given by
$$
\int_{X^\GG_k(\Theta_{\cL_k,h_{V,F,k,\varepsilon}},q)\ssm\pi_k^{-1}(\Sigma)}
\Theta_{\cL_k,h_{V,F,k,\varepsilon}}^{\raise2.5pt\hbox{$\scriptstyle n+kr-1$}}=
{(\log k)^n\over n!\,(k!)^r}\bigg(
\int_{X(\eta,q)\ssm\Sigma}\bOne_{\eta,q}\eta^n+O((\log k)^{-1})\bigg)
\leqno(5.6\,{\rm d})
$$
for all $q=0,1,\ldots,n$, and the error term $O((\log k)^{-1})$ can be 
bounded explicitly in terms of $\Theta_{V,h}$ and $\Theta_{F,h_F}$.
Moreover, the left hand side is identically zero for $q>n$.
\endclaim

\proof We refer to [Dem11] which contains all details in the smooth case,
i.e.\ when $V\subset T_X$ is a subbundle, and when $h_V$, $h_F$ are smooth.
The main argument is that the $(1,1)$ curvature form of $\cO_{X_k^\GG}(1)$
can be expressed explicitly, modulo small error terms. It is convenient
to use polar coordinates
$$
\xi_s=\varepsilon_s^{-1/2p}x_s^{1/p}u_s,\leqno(5.7)
$$
where $u_s\in SV$ is in the unit sphere bundle of $V$ ($|u_s|_h=1)$ and
$x_s\geq 0$. Then at any point $(z,\xi)\in X_k^\GG$, a straightforward
calculation shows that the curvature of $\cO_{X_k^\GG}(1)$ is given by
$$
\Theta_{\cO_{X_k^\GG}(1),h_{V,k,\varepsilon}}(z,\xi)=
\omega_{p,\FS}(\xi)-{i\over 2\pi}
\sum_{1\leq s\leq k}{x_s\over s}\sum_{i,j,\alpha,\beta}
c_{ij\alpha\beta}(z)\,u_{s\alpha}\overline{u_{s\beta}}\,dz_i\wedge d\overline z_j
+O(\varepsilon)\leqno(5.8)
$$
where
$$
\omega_{p,\FS}(\xi)={1\over p}{i\over 2\pi}\ddbar\log
\sum_{1\leq s\leq k}|\xi_s|^{2p/s}
$$
is the weighted Fubini-Study metric on the fibers of $X_k^\GG\to X$.
As $k\to+\infty$, the double summation of (5.8) can be seen as
a Monte-Carlo evaluation of the curvature tensor
$u\mapsto \langle\Theta_{V,h_V}u,u\rangle$ on
the sphere bundle $SV$. It thus exhibits a probabilistic
convergence, and is on average equivalent to a quantity essentially
independent of $\xi$, namely
$$
-{1\over r}\sum_{1\leq s\leq k}{x_s\over s}\Tr_{\End(V)}\Theta_{V,h_V}
$$
proportional to $\Theta_{\det(V^*),\det(h^*)}$ (the integral of
a quadratic form on a sphere is just its trace!). Then (5.6$\,$d) follows by
analyzing the error terms and computing fiber integrals (the latter
depend only on the weighted Fubini-Study metric and can be easily evaluated).
One just needs to add $\Theta_{F,h_F}$ to estimate the Morse integral of
$\cL_k$. The corresponding integrals vanish for $q>n$ because $\cL_k$
is semi-positive (and generically strictly positive) along the
fibers of $\pi_k:X_k^\GG\to X$, and the only negative eigenvalues are
those coming from $\eta(z)$ in the $dz_i$'s. One important point is that
the rescaling factor $\varepsilon_s$ in (5.8) makes the metric
(5.6$\,$b) essentially independent of the connection~$\nabla_0$,
up to errors $O(\varepsilon)$ (cf.\ Lemma~2.12 in [Dem11]). Therefore,
one can use a local holomorphic connection $\nabla$ of $V$ (say, a trivial
flat one) to perform calculations -- this is helpful to ensure that
the $\xi_s$ are actually holomorphic coordinates. In the singular case,
we use a smooth modification $\mu:\widetilde X\to X$ that takes
$V$ to a locally free sheaf $\mu^*V$ and converts $h_V$, $h_F$ into
metrics with $\bQ$-divisorial singularities on $\smash{\widetilde X}$,
in such a way that we are reduced to the smooth case, after removing
the corresponding divisors and computing the Morse integrals in the
complement of~$\Sigma$. However,
it must be observed that we still need in (5.6$\,$b) to use a {\it smooth}
(or local holomorphic) connection $\nabla_0$ on $T_X$, and not a
holomorphic connection $\nabla_V$ on $V$, because such a connection would
``explode'' near the singularities. Far away from the singularities,
there is uniform convergence to what would be obtained from a holomorphic
connection $\nabla_V$ on~$V$. Near $\Sing(V)$, the argument is that
the curvature of $\Theta_{V,h_X}$ is controlled by the
second fundamental form of~$V$ (the relevant term being negative), but
the total volume of any exterior power is convergent and therefore goes
to $0$ if we pick a sufficiently small neighborhood of the singular set
$\Sing(V)$. This would not work if we had chosen a holomorphic
connection on~$V$ near the singularities! Also notice that the singular
factor $e^\varphi$ of $h=e^\varphi(h_X)_{|V}$ does not interfere as it
is purely scalar and ``diagonal''.\qed

Our general Morse inequalities 4.7 then yield directly

\claim 5.9. Morse inequalities for tautological sheaves|
Let $(X,V)$ be a directed variety, where $X$~is a nonsingular compact
complex manifold, let $F$ be a $\bQ$-line bundle and $\cG$
a coherent sheaf of rank~$\rho$ on~$X$. Assume that $V$ and $F$ are equipped
with metrics $h_V$ and $h_F$ with $\bQ$-analytic singularities
and let $\eta=\Theta_{\det(V^*),\det(h^*)}+\Theta_{F,h_F}$,
$X'=X\ssm(\Sing(h_V)\cup\Sing(h_F))$. Then 
for $m\gg k\gg 1$ we have estimates
$$
\eqalign{
\rho\,{m^{n+kr-1}\over (n+kr-1)!}
{(\log k)^n\over n!\,(k!)^r}\Bigg(&\sum_{j=q-1,q,q+1}
\int_{X'(\eta,j)}(-1)^{q-j}\eta^n-O((\log k)^{-1})\Bigg),\cr
\leq h^q\bigg(X_k^\GG,{}^b&\cO_{X_k^\GG}(m)_{h_{k,\varepsilon}^{-m}}\otimes
\pi_k^*\Big({}^b\cO\Big({m\over kr}\Big(1+{1\over 2}+\ldots+{1\over k}
\Big)F\Big)_{\raise2pt\hbox{$h_F$}}\otimes\cG\Big) \bigg)\cr
&\kern30pt{}\leq \rho\,{m^{n+kr-1}\over (n+kr-1)!}
{(\log k)^n\over n!\,(k!)^r}\Bigg(
\int_{X'(\eta,q)}(-1)^q\eta^n+O((\log k)^{-1})\Bigg).\cr}
$$
\endclaim

\proof Apply 4.7 on $Z=X_k^\GG$ to the $m$-th power of the rank $1$ sheaf
$$
\cL_k=\cO_{X_k^\GG}(1)\otimes
\pi_k^*\bigg(\cO_X\bigg({1\over kr}\Big(1+{1\over 2}+\ldots+{1\over k}
\Big)F\bigg)\bigg)
$$
equipped with the metric induced by $(h_{V,k,\varepsilon})^{-1}$ and $h_F$,
and evaluate the resulting cohomology of the sheaf of bounded sections of
$\cL_k^m\otimes\pi_k^*\cG$.\qed

Notice that there is in fact no loss of
generality to take $h_V=(h_X)_{|V}$ where $h_X$ is a smooth metric on~$T_X$,
since any additional weight $e^{\varphi}$ with $\bQ$-analytic singularities
can in fact be moved to the factor~$F$. In this case, we simply denote by
${}^b\cO_{X_k^\GG}(m)$ the sheaf of bounded sections and do not specify
the metric on~$V$. By taking the direct image via $\pi_k:X_k^\GG\to X$
and applying the Leray spectral sequence
we obtain the same bounds as above for the cohomology groups
$$
H^q\bigg(X,{}^b\cO(E_{k,m}^\GG V^*)\otimes
{}^b\cO\Big({m\over kr}\Big(1+{1\over 2}+\ldots+{1\over k}
\Big)F\Big)_{\raise2pt\hbox{$h_F$}}\otimes\cG\bigg),
$$
putting now $\eta=\Theta_{\det(V^*),\det(h_X^*)}+\Theta_{F,h_F}$.
By Serre duality, we infer similar bounds for the ``dual''
cohomology groups
$$
H^{n-q}\bigg(X,\big({}^b\cO(E_{k,m}^\GG V^*)\big)^*\otimes
{}^b\cO\Big(-{m\over kr}\Big(1+{1\over 2}+\ldots+{1\over k}
\Big)F\Big)_{\raise2pt\hbox{$h_F$}}\otimes\cG'\bigg)
$$
with $\cG'=K_X\otimes\cG^*$, but the duality holds only
when $\cG$ and ${}^b\cO(E_{k,m}^\GG V^*)$ are locally free,
which is the case if $V$ is nonsingular. In fact, there is no
change for the dominant term of the bounds if $\cG$ is not locally free,
because we can replace $\cG$ by $K_X\otimes\cG^*$ and then
$\cG'$ equals $\cG^{**}$, which coincides with $\cG$ outside
of codimension~$1$ (cf.\ the proof of Lemma~4.4), and anyway the
dominant term depends only on the generic rank
$\rho=\rank\cG$. When $V$ is singular, we have by definition an injection
$J^kV\hookrightarrow J^kT_X$. It~yields a restriction morphism
$$
\cO(E_{k,m}^\GG T_X^*)={}^b\cO(E_{k,m}^\GG T_X^*)\to{}^b\cO(E_{k,m}^\GG V^*),
$$
and by duality an injection
$$
\big({}^b\cO(E_{k,m}^\GG V^*)\big)^*\hookrightarrow\cO(E_{k,m}^\GG T_X^*)^*
\leqno(5.10)
$$
where the right hand side is a vector bundle. In order to deal properly with
the duality, one way is use a modification $\mu:\widetilde X\to X$
such that $\mu^*({}^b\cO(E_{k,m}^\GG T_X^*))$ is locally free modulo torsion.
This modification can be chosen independent of $m$ because at any point $x$,
$\bigoplus_m{}^b\cO_{X,x}(E_{k,m}^\GG V^*)$ is a finitely generated 
graded algebra over $\cO_{X,x}$ (equal to the integral closure of the image
of the finitely generated algebra $\bigoplus_m\cO_{X,x}(E_{k,m}^\GG T_X^*)$
in its field of quotients). Once this is done, one can instead pull-back
to $\smash{\widetilde X}$ and apply Serre duality on $\smash{\widetilde X}$.
Again, the dominant term is unchanged as the replacement of $K_X$ by
$K_{\smash{\widetilde X}}$
leaves it unaffected. This is especially interesting for $q=n$ since
we then obtain an estimate of holomorphic sections of
$\big(\mu^*\,{}^b\cO(E_{k,m}^\GG V^*)\big)^*$ on $\smash{\widetilde X}$,
and therefore of
$\smash{\mu_*\big(\big(\mu^*\,{}^b\cO(E_{k,m}^\GG V^*)\big)^*\big)}$
on~$X$. If we reinterpret these operations in terms of metrics and bounded
sections, we obtain precisely the sheaf of bounded sections
$\smash{{}^b\cO\big((E_{k,m}^\GG V^*)^*\big)}$, defined as the space of sections
of $(E_{k,m}^\GG V^*)^*$ on $X\ssm\Sing(V)$ that are bounded with respect
to a smooth metric on $\smash{(E_{k,m}^\GG T_X^*)^*}$, and we see that we have
in fact
$$
{}^b\cO\big((E_{k,m}^\GG V^*)^*\big)=
\mu_*\big(\big(\mu^*\,{}^b\cO(E_{k,m}^\GG V^*)\big)^*\big)=
\big({}^b\cO(E_{k,m}^\GG V^*)\big)^*.
$$
From these considerations we infer

\claim 5.11. Morse inequalities for holomorphic jet differentials and
their duals|
Let $(X,V)$ be a directed variety, where $X$ is a nonsingular compact
complex manifold, let $F$ be a $\bQ$-line bundle and $\cG$
a coherent sheaf of rank $\rho$ on~$X$. Assume that $F$ is equipped
with a metric $h_F$ with $\bQ$-analytic singularities and let $h_X$ be
a smooth metric on $T_X$. For $m\gg k\gg 1$ we have the following estimates.
\plainitem{\rm (a)}Let $\eta=\Theta_{\det(V^*),\det(h_X^*)}+\Theta_{F,h_F}$
and $X'=X\ssm(\Sing(V)\cup\Sing(h_F))$. Then
$$
\eqalign{
\rho\,{m^{n+kr-1}\over (n+kr-1)!}
{(\log k)^n\over n!\,(k!)^r}\Bigg(&\sum_{j=q-1,q,q+1}
\int_{X'(\eta,j)}(-1)^{q-j}\eta^n-O((\log k)^{-1})\Bigg),\cr
\leq h^q\bigg(X,{}^b&\cO(E^\GG_{k,m}V^*)\otimes
{}^b\cO\Big({m\over kr}\Big(1+{1\over 2}+\ldots+{1\over k}
\Big)F\Big)_{\raise2pt\hbox{$h_F$}}\otimes\cG\bigg)\cr
&\kern30pt{}\leq \rho\,{m^{n+kr-1}\over (n+kr-1)!}
{(\log k)^n\over n!\,(k!)^r}\Bigg(
\int_{X'(\eta,q)}(-1)^q\eta^n+O((\log k)^{-1})\Bigg).\cr}
$$
\plainitem{\rm (b)}Let $\eta^*=\Theta_{\det(V),\det(h_X)}+\Theta_{F,h_F}$ and
$X'=X\ssm(\Sing(V)\cup\Sing(h_F))$. Then
$$
\eqalign{
\rho\,{m^{n+kr-1}\over (n+kr-1)!}
{(\log k)^n\over n!\,(k!)^r}\Bigg(&\sum_{j=q-1,q,q+1}
\int_{X'(\eta^*,j)}(-1)^{q-j}(\eta^*)^n-O((\log k)^{-1})\Bigg),\cr
\leq h^q\bigg(X,{}^b&\cO\big((E^\GG_{k,m}V^*)^*\big)\otimes
{}^b\cO\Big({m\over kr}\Big(1+{1\over 2}+\ldots+{1\over k}
\Big)F\Big)_{\raise2pt\hbox{$h_F$}}\otimes\cG\bigg)\cr
&\kern30pt{}\leq \rho\,{m^{n+kr-1}\over (n+kr-1)!}
{(\log k)^n\over n!\,(k!)^r}\Bigg(
\int_{X'(\eta^*,q)}(-1)^q(\eta^*)^n+O((\log k)^{-1})\Bigg).\cr}
$$
\endclaim

\proof (a) is a consequence of our direct image argument. part (b) follows by
duality if we observe that $X(\eta,n-q)=X(q,-\eta)$ and change $F$
into $-F$; this has the effect of replacing $-\eta$ by $\eta^*$ and
$(-1)^{n-q}\eta^n$ by $(-1)^q(\eta^*)^n$.\qed

\claim 5.12. Definition|Let $\cL$ be a rank one  sheaf equipped
with a hermitian metric $h_\cL$ with $\bQ$-analytic singularities. We say 
for short that ${}^b\cL$ is big if there exists $m\in\bN^*$ and a 
log resolution 
$\smash{\mu:\widetilde X\to X}$ of the sheaf ${}^b\cL^m={}^b\cL^m_{h_\cL}$ 
of bounded sections such that $\mu^*({}^b\cL^m)$ is a big line bundle.
When taking the product $\cL\otimes\cF$ with an invertible sheaf and speaking
of the bigness of $\cL\otimes\cF$, we agree implicitly to take the tensor 
product $h_\cL\otimes h_\cF$ with a smooth metric on $\cF$, so that
${}^b(\cL\otimes\cF)^m=({}^b\cL^m)\otimes\cF^m$.
\endclaim

With this terminology in mind, we can now state an important application of 
our Morse esti\-mates, in relation with positivity properties of the
canonical or anticanonical sheaf of a directed variety.

\claim 5.13. Corollary|Let $X$ be a projective $n$-dimensional manifold,
$(X,V)$ a directed structure, and $F$ an invertible sheaf on $X$.
We consider here bounded sections with respect to a smooth hermitian
metric $h_X$ on~$T_X$.
\plainitem{\rm(a)} If ${}^bK_V\otimes F$ is big, there are many sections
in
$$
H^0\bigg(X,{}^b\cO\big(E^\GG_{k,m}V^*\big)\otimes
{}^b\cO\Big({m\over kr}\Big(1+{1\over 2}+\ldots+{1\over k}
\Big)F\Big)\bigg)\quad\hbox{for $m\gg k\gg 1$}.
$$
\plainitem{\rm(b)} If ${}^b\cO(\det(V))\otimes F$ is big, there are
many sections in
$$
H^0\bigg(X,{}^b\cO\big((E^\GG_{k,m}V^*)^*\big)\otimes
{}^b\cO\Big({m\over kr}\Big(1+{1\over 2}+\ldots+{1\over k}
\Big)F\Big)\bigg)\quad\hbox{for $m\gg k\gg 1$}.
$$
The asymptotic growth is of the form $c\,m^{n+kr-1}( \log k)^n/
((n+kr-1)!\,(k!)^r)$ where $r=\rank V$ and the constant $c>0$ depends only
on $X$, $V$ and $F$.
\endclaim

\proof In case (a), by [Dem90], we can find a singular metric $h_F$
on $F$ such that
$$
\eta=\Theta_{K_V,\det(h_X)^*)}+\Theta_{F,h_F}
$$
is a K\"ahler current, i.e.\ $\eta\geq c\omega$ for some
K\"ahler metric $\omega$ on $X$ and some small constant~$c>0$.
In fact, we can work with the invertible sheaf $\mu^*({}^bK_V\otimes F)$
on $\smash{\widetilde X}$ and observe that the push forward
of a K\"ahler current on $\smash{\widetilde X}$ is a K\"ahler current on $X$.
Then all chambers $X'(\eta,q)$ are empty except $X'(\eta,0)$
which yields a strictly positive Morse integral, whence the
result. Part (b) is entirely similar, using $\eta^*$ instead
of $\eta$. Notice that this does not necessarily imply that the
higher cohomology groups vanish, only that
$h^q_{m,k}=O((\log k)^{-1})h^0_{m,k}$ as~$m\gg k\gg 1$.\qed

\claim 5.14. Remark|{\rm Since $\smash{\cO(E^\GG_{k,m}V^*)}$ is the sheaf
of homogeneous polynomials of degree $m$ on $J^kV$, its
dual $\smash{\cO((E^\GG_{k,m}V^*)^*)}$ can be thought of as a sheaf of
differential operators of degree $m$ on $J^kV$. Hence, our result
(5.13$\,$b) can be seen as an extension of the results of
Siu [Siu04], P\u{a}un [Pau08], Merker [Mer09] on the existence of
slanted vector fields on jet bundles.}
\endclaim \bigskip

\section*{References}
\vskip5pt

\begingroup
\fontsize{10pt}{12pt}\selectfont

\Bibitem[Blo26a]&Bloch, A.:& Sur les syst\`emes de fonctions uniformes
satisfaisant \`a l'\'equation d'une vari\'et\'e alg\'ebrique dont
l'irr\'egularit\'e d\'epasse la dimension.& J.\ de Math., {\bf 5}
(1926), 19--66&

\Bibitem[Blo26b]&Bloch, A.:& Sur les syst\`emes de fonctions holomorphes
\`a vari\'et\'es lin\'eaires lacunaires.& Ann.\ Ecole Normale, {\bf 43}
(1926), 309--362&

\Bibitem[Bon93]&Bonavero, L.:& In\'egalit\'es de Morse holomorphes
singuli\`eres.& C.~R.\ Acad.\ Sci.\ Paris S\'er.~I Math.\ {\bf 317} (1993)
1163--1166&

\Bibitem[Dem85]&Demailly, J.-P.:& Champs magn\'etiques et in\'egalit\'es de
Morse pour la $d''$-coho\-mo\-logie.& Ann.\ Inst.\ Fourier
(Grenoble) {\bf 35} (1985) 189--229&

\Bibitem[Dem90]&Demailly, J.-P.:& Singular hermitian metrics on positive
line bundles.& Proceedings of the Bayreuth conference ``Complex algebraic
varieties'', April~2--6, 1990, edited by K.~Hulek, T.~Peternell,
M.~Schneider, F.~Schreyer, Lecture Notes in Math.\ ${\rm n}^\circ\,$1507,
Springer-Verlag (1992), 87--104&

\Bibitem[Dem11]&Demailly, J.-P.:& 
Holomorphic Morse Inequalities and the Green-Griffiths-Lang Conjecture.&
Pure and Applied Math.\ Quarterly {\bf 7} (2011), 1165--1208&

\Bibitem[Dem12]&Demailly, J.-P.:& 
Hyperbolic algebraic varieties and holomorphic differential equations.&
expanded version of the lectures given at the annual meeting of VIASM,
Acta Math.\ Vietnam.\ {\bf 37} (2012), 441-–512&

\Bibitem[Dem15]&Demailly, J.-P.:&
Towards the Green-Griffiths-Lang conjecture,&
Conference ``Analysis and Geo\-metry'', Tunis,
March 2014, in honor of Mohammed Salah Baouendi, 
ed.\ by A.~Baklouti, A.\ El Kacimi, S.~Kallel, N.~Mir, Springer (2015),
141--159&

\Bibitem[Dem18]&Demailly, J.-P.:&
Recent results on the Kobayashi and Green-Griffiths-Lang conjectures,&
arXiv: Math.AG/1801.04765, Contribution to the 16th Takagi lectures
in celebration of the 100th anniversary of K.Kodaira's birth,
November 2015, to appear in the Japanese Journal of Mathe\-matics&

\Bibitem[GrGr79]&Green, M., Griffiths, P.:& Two applications of algebraic
geometry to entire holomorphic mappings.& The Chern Symposium 1979,
Proc.\ Internal.\ Sympos.\ Berkeley, CA, 1979, Springer-Verlag, New York
(1980), 41--74&

\Bibitem[Lang86]&Lang, S.:& Hyperbolic and Diophantine analysis.&
Bull.\ Amer.\ Math.\ Soc.\ {\bf 14} (1986), 159--205&

\Bibitem[Mer09]&Merker, J.:& Low pole order frames on vertical jets of the 
universal hypersurface& Ann.\ Inst.\ Fourier (Grenoble), {\bf 59} (2009),
1077--1104&

\Bibitem[Pau08]&P\u{a}un, M.:& Vector fields on the total space of 
hypersurfaces in the projective space and hyperbolicity.&
Math.\ Ann.\ {\bf 340} (2008) 875--892&

\Bibitem[Siu04]&Siu, Y.T.:& Hyperbolicity in complex geometry.& The 
legacy of Niels Henrik Abel, Springer, Berlin (2004) 543--566&

\endgroup

\vskip15pt
\parindent=0cm
(version of June 25, 2018, printed on \today, \timeofday)
\vskip15pt

Jean-Pierre Demailly\\
Universit\'e de Grenoble-Alpes, Institut Fourier (Math\'ematiques)\\
UMR 5582 du C.N.R.S., 100 rue des Maths, 38610 Gi\`eres, France\\
{\em e-mail}$\,$: jean-pierre.demailly@univ-grenoble-alpes.fr
\medskip


Mohammad Reza Rahmati\\
Institute of Algebraic Geometry, Gottfried 
Wilhelm Leibniz Universit\"at Hannover\\
Welfengarten 1, 30167 Hannover, Germany\\
{\em e-mail}$\,$: mrahmati.mrr@gmail.com

\end{document}